\documentclass{amsart}

\usepackage{amsthm}
\usepackage{amssymb}
\usepackage{upref}
\usepackage{amscd}

\numberwithin{equation}{section}

 \newtheorem{theorem}[equation]{Theorem}
 
 \newtheorem{proposition}[equation]{Proposition}
 \newtheorem{corollary}[equation]{Corollary}
 
 \newtheorem{question}[equation]{Question}
 \theoremstyle{definition}
 \newtheorem{definition}[equation]{Definition}

 \DeclareMathOperator{\Hom}{Hom}
 \DeclareMathOperator{\Ext}{Ext}
 \DeclareMathOperator{\End}{End}
 \DeclareMathOperator{\Tor}{Tor}
 \DeclareMathOperator{\colim}{colim}
 \DeclareMathOperator{\ghdim}{gh.dim.}
 
 \DeclareMathOperator{\wdim}{w. dim.}
 \DeclareMathOperator{\rgldim}{r. gl. dim.}
 \DeclareMathOperator{\pgldim}{pure \ gl. dim.}
 \DeclareMathOperator{\gldim}{gl. dim.}
 \DeclareMathOperator{\pdim}{proj. dim.}
 \DeclareMathOperator{\idim}{inj. dim.}
 \DeclareMathOperator{\fdim}{flat \  dim.}
 \DeclareMathOperator{\cfdim}{con.\ flat \  dim.}
 \DeclareMathOperator{\phdim}{phan. dim.}
 \DeclareMathOperator{\Sq}{Sq}
 \DeclareMathOperator{\depth}{depth}

\newcommand{\Q}{\mathbb{Q}}
\newcommand{\ideal}[1]{\mathfrak{#1}}
\newcommand{\cat}[1]{\mathcal{#1}}
\newcommand{\Z}{\mathbb{Z}}

\newcommand{\mathcolon}{\colon\,}

\newcommand{\uc}{\textup{:}}
\newcommand{\ulp}{\textup{(}}
\newcommand{\urp}{\textup{)}}

\begin{document}

\title{Homological dimensions of ring spectra} 

\date{\today}

\author{Mark Hovey}
\address{Department of Mathematics \\ Wesleyan University
\\ Middletown, CT 06459}
\email{hovey@member.ams.org}

\author{Keir Lockridge}
\address{Department of Mathematics \\ Wake Forest University \\
Winston-Salem, NC 27109}
\email{lockrikh@wfu.edu}


\begin{abstract}
We define homological dimensions for $S$-algebras, the generalized
rings that arise in algebraic topology.  We compute the homological
dimensions of a number of examples, and establish some basic
properties.  The most difficult computation is the global dimension of
real $K$-theory $KO$ and its connective version $ko$ at the prime
$2$.  We show that the global dimension of $KO$ is $1$, $2$, or $3$,
and the global dimension of $ko$ is $4$ or $5$.  
\end{abstract}

\maketitle

\section*{Introduction}

The authors have been engaged in trying to develop homological
dimensions for the ring objects that arise in algebraic topology.
These are cohomology theories with some kind of cup product that is
associative up to infinitely coherent homotopy.  These are commonly
called $S$-algebras, among other names, and a standard reference
is~\cite{elmendorf-kriz-mandell-may}.

Such $S$-algebras are analogous to rings, but they have no elements,
and they do not have abelian module categories.  What they do have is
a triangulated derived category that generalizes the derived category
of an ordinary associative ring $R$.  Indeed, in this case, there is
an Eilenberg-MacLane spectrum $S$-algebra $HR$, and $\cat{D} (HR)$ is
naturally equivalent to $\cat{D} (R)$ as a triangulated category.

Thus, to develop the ring theory of $S$-algebras $E$, we need to work
with their derived categories $\cat{D} (E)$.  In previous
papers~\cite{hovey-lockridge-semisimple,hovey-lockridge-weak}, the
authors have developed definitions of right global dimension $\rgldim
E$ and ghost dimension $\ghdim E$ for $S$-algebras $E$ and proved that
these generalize the usual notions of right global dimension and weak
dimension for rings.  We have also studied $S$-algebras for which
these dimensions are $0$.  

The object of this paper is to study the homological dimensions of the
$S$-algebras that arise in nature.  After an initial section in which
we recall the basics of our theory of global dimension and ghost
dimension of $S$-algebras, we discuss examples in
Section~\ref{sec-examples}.  The sphere, $MU$, and $BP$ all have
infinite global and ghost dimensions, whereas $E_{n}$ has ghost and
global dimension $n$.  For a commutative $S$-algebra $E$ such that $E_{*}$
is Noetherian and has finite global dimension, such as $E_{n}$ or $K$,
we have 
\[
\gldim E=\ghdim E =\gldim E_{*}.
\]
This fails when $\gldim E_{*}=\infty$, however.  We show that, at the
prime $2$, $\gldim KO$ is either $1$, $2$, or $3$, and $\gldim ko$ is
either $4$ or $5$, whereas both rings $KO_{*}$ and $ko_{*}$ have
infinite global dimension.  We also show that $\gldim tmf$ is finite
at the prime $3$.

In Section~\ref{sec-properties} we prove some basic properties of
global and ghost dimension.  The most important point is that these
dimensions are not Morita invariant.  In fact, it is already true that
the global and weak dimensions of ordinary rings are not invariant
under derived equivalences, and we give an example communicated to us
by Lidia Angeleri H\"{u}gel.  We have been unable to generalize many
of the basic properties of the global dimension of Noetherian rings,
however.  This might be because we have no intrinsic definition of $E$
being a Noetherian $S$-algebra; we just assume $E_{*}$ is right
Noetherian.  So, for example, we do not know whether $\rgldim E=\ghdim
E$ when $E_{*}$ is right Noetherian.  We end the paper with a brief
section on $S$-algebras $E$ with global dimension $1$.  We had
originally thought this would mean $E_{*}$ would have to be
$1$-Gorenstein, in analogy with the fact that $\gldim E=0$ implies
$E_{*}$ is quasi-Frobenius, and in fact this is true, but only with
additional assumptions on $E_{*}$.

The authors would like to thank Lidia Angeleri H\"{u}gel for the
example mentioned above, and Ben Wieland for pushing us to determine
the global dimension of $KO$.

\section{Homological dimensions}

The object of this section is to define our various homological
dimensions of ring spectra, and to prove basic relations between
them.  Recall that $E$ is an $S$-algebra, and $\cat{D} (E)$ denotes
the derived category of right $E$-modules.  This is a compactly
generated triangulated category, and when $E$ is the Eilenberg-MacLane
spectrum $HR$ of an ordinary ring $R$, then $\cat{D} (E)$ is
equivalent to the usual unbounded derived category of $R$.  

We begin by defining projective, injective, and flat objects of
$\cat{D} (E)$.  

\begin{definition}
An object $X\in \cat{D} (E)$ is said to be \textbf{projective}
(resp. \textbf{injective}, resp. \textbf{flat}) if $X_{*}$ is a
projective (resp. injective, resp. flat) $E_{*}$-module.  A map
$f\mathcolon X\xrightarrow{}Y$ in $\cat{D} (E)$ is said to be
\textbf{ghost} if $f_{*}=0$, and $f$ is said to be \textbf{phantom}
if $\cat{D} (E) (A,f)=0$ for all compact $A\in \cat{D} (E)$.  
\end{definition}
\label{defn-projective}

The basic properties of these objects and maps are summed up in the
following proposition.

\begin{proposition}\label{lem-projective}
Suppose $E$ is an $S$-algebra.  
\begin{enumerate}
\item If $M$ is a projective or injective $E_{*}$-module, then
there is an $X\in \cat{D} (E)$ \ulp necessarily projective or
injective\urp  with $X_{*}\cong M$.  
\item If $M$ is a flat $E_{*}$-module, and $\cat{D} (E)$ is a Brown
category \ulp see~\cite[Section~4.1]{hovey-axiomatic}\urp, for example
if $E_{*}$ is countable, then there is an $X\in \cat{D} (E)$ \ulp
necessarily flat\urp  such that $X_{*}\cong M$.
\item $X$ is a projective object of $\cat{D} (E)$ if and only if the
natural map 
\[
\cat{D} (E) (X,Y)\xrightarrow{}\Hom_{E_{*}}
(X_{*},Y_{*})
\]
is an isomorphism for all $Y\in \cat{D} (E)$.  This is true if and
only if every ghost with domain $X$ is null.
\item $X$ is an injective object of $\cat{D} (E)$ if and only if the
natural map 
\[
\cat{D} (E) (Y,X)\xrightarrow{}\Hom_{E_{*}}
(Y_{*},X_{*})
\]
is an isomorphism for all $Y\in \cat{D} (E)$.  This is true if and
only if every ghost with codomain $X$ is null.
\item $X$ is a flat object of $\cat{D} (E)$ if and only the natural
map 
\[
X_{*}\otimes_{E_{*}}Y_{*}\xrightarrow{}\pi_{*} (X\wedge_{E}Y)
\]
is an isomorphism for all left $E$-modules $Y$.  This is true if and
only if every ghost with domain $X$ is phantom.
\end{enumerate}
\end{proposition}

\begin{proof}
The first part is well-known.  In fact, it is proved
in~\cite[Proposition~A.4]{benson-krause-schwede} that every
$E_{*}$-module of projective or injective dimension $\leq 2$ is
realizable.  We prove part~(3) next.  If $X$ is projective, then the
universal coefficient spectral sequence
of~\cite[Theorem~IV.4.1]{elmendorf-kriz-mandell-may} implies that 
\[
\cat{D} (E) (X,Y)\cong \Hom_{E_{*}} (X_{*},Y_{*}).
\]
This in turn implies that there are no nontrivial ghosts with domain
$X$.  Now, if there are no nontrivial ghosts with domain $X$,
construct a projective module $P_{*}$ and an epimorphism
$P_{*}\xrightarrow{}X_{*}$.  This is then realizable, as we have just
seen, by a map $P\xrightarrow{}X$, whose cofiber $X\xrightarrow{}Y$ is
a ghost.  This map is thus null, so $X$ is a retract of $P$, and hence
projective.  The proof of part~(4) is similar. Part~(5) is proved
in~\cite{hovey-lockridge-weak}.

Turning to part~(2), suppose $M$ is a flat $E_{*}$-module.  Then
$M_{*}$ is a directed colimit of finitely generated projective
$E_{*}$-modules $Q_{i}$.  We can then realize this system by compact
projective $E$-modules $P_{i}$, with $\pi_{*}P_{i}\cong Q_{i}$, by
part~(3).  If $\cat{D} (E)$ is a Brown category, we can then take the
minimal weak colimit~\cite[Section~4.2]{hovey-axiomatic} of the
diagram of the $P_{i}$ to obtain an $X$ with $X_{*}\cong \colim
Q_{i}\cong M$.
\end{proof}

Once we have projective, injective, and flat objects, we should then
be able to define the projective dimension of an arbitrary object.
This is a little more complicated than in the abelian setting, but has
in fact already been worked out by Christensen in~\cite{christensen}.  

\begin{definition}\label{defn-dimension}
Let $E$ be an $S$-algebra, and $X$ an object of $\cat{D} (E)$.  We
define the \textbf{projective dimension} (resp. \textbf{constructible
flat dimension}) of $X$, written $\pdim X$ (resp. $\cfdim X$),
inductively as follows.  We have $\pdim X=0$ (resp. $\cfdim X=0$) if
and only if $X$ is projective (resp. flat).  Then we define $\pdim
X\leq n+1$ (resp. $\cfdim X\leq n+1$) if and only if there is an exact
triangle
\[
Y \xrightarrow{} P \xrightarrow{} \widetilde{X}\xrightarrow{}\Sigma Y
\]
where $P$ is projective (resp. flat), $\pdim Y\leq n$ (resp. $\cfdim
X\leq n$), and $X$ is a retract of $\widetilde{X}$. We define the
\textbf{flat dimension} of $X$, written $\fdim X$, to be the smallest
integer $n$ for which any composite of $n+1$ ghosts with domain $X$ is
phantom, or $\infty$ if there is no such $n$.  The reader may well
wonder why we have two notions of flat dimension.  The answer is that
we have been unable to prove they are equivalent, and the
constructible flat dimension is the more obvious one, but the flat
dimension is the more useful one.  See Proposition~\ref{prop-proj-dim}
below and~\cite{hovey-lockridge-weak} for details.

We can define injective dimension similarly to how we defined
projective dimension, but it would be more usual to define the
\textbf{injective dimension} of $X$, written $\idim X$, inductively as
follows.  We define $\idim X=0$ if and only if $X$ is injective, and
$\idim X\leq n+1$ if and only if there is an exact triangle
\[
\Sigma^{-1}Y \xrightarrow{}\widetilde{X} \xrightarrow{} I \xrightarrow{} Y
\]
where $I$ is injective, $\idim Y\leq n$, and $X$ is a retract of
$\widetilde{X}$.  
\end{definition}

The major difference between this definition and the definition of
the analogous dimensions in abelian categories is the fact that $X$
itself need not appear in an exact triangle with things of smaller
dimension, but must only be a retract of such a thing.  Without this
condition, Proposition~\ref{prop-proj-dim} below would be false.  

\begin{proposition}\label{prop-proj-dim}
Let $E$ be an $S$-algebra, and let $X$ be an object of $\cat{D} (E)$.  
\begin{enumerate}
\item $\pdim X\leq n$ if and only if every composite of $n+1$ ghosts
$f_{n+1}\circ f_{n}\circ \dotsb \circ f_{1}$ is null, where the domain
of $f_{1}$ is $X$.  This is true if and only if in the universal
coefficient spectral sequence
\[
E_{2}^{s,t} = \Ext_{E_{*}}^{s,t} (X_{*},Y_{*}) \Rightarrow \cat{D} (E)
(X,Y)_{t-s}
\]
we have $E_{\infty}^{s,*}=0$ for all $s>n$ and all objects $Y$ of
$\cat{D} (E)$.   
\item $\idim X\leq n$ if and only if every composite of $n+1$ ghosts
$f_{n+1}\circ f_{n}\circ \dotsb \circ f_{1}$ is null, where the
codomain of $f_{n+1}$ is $X$.  This is true if and only if in the
universal coefficient spectral sequence
\[
E_{2}^{s,t} = \Ext_{E_{*}}^{s,t} (Y_{*},X_{*}) \Rightarrow \cat{D} (E)
(Y,X)_{t-s}
\]
we have $E_{\infty}^{s,*}=0$ for all $s>n$ and all objects $Y$ of
$\cat{D} (E)$.
\item $\fdim X\leq \cfdim X$, and $\fdim X \leq n$ if and only if the
following equivalent properties hold\uc
\begin{enumerate}
\item In the universal coefficient spectral sequence
\[
E^{2}_{s,t} = \Tor^{E_{*}}_{s,t} (X_{*},Y_{*}) \Rightarrow \pi_{t-s}
(X\wedge_{E}Y)
\]
we have $E^{\infty}_{s,*}=0$ for all $s>n$ and all objects $Y$ of
$\cat{D} (E^{\textup{op}})$.
\item There is an exact triangle 
\[
A \xrightarrow{}X \xrightarrow{g} W\xrightarrow{} \Sigma A
\]
in which $\pdim A\leq n$ and $g$ is phantom.  
\item Every map $F\xrightarrow{}X$ from a compact $E$-module $F$
factors through a compact $B$ with $\pdim B\leq n$.  
\end{enumerate}
\end{enumerate}
\end{proposition}

\begin{proof}
For part~(1), the first part is contained
in~\cite[Theorem~3.5]{christensen}.  The second part follows from the
construction and naturality of the universal coefficient spectral
sequence~\cite[Section~IV.5]{elmendorf-kriz-mandell-may}.  Indeed,
suppose $E_{\infty}^{s,*}=0$ for all $s>n$ and all $Y$, and suppose 
\[
X\xrightarrow{f_{1}} X_{1} \xrightarrow{f_{2}}\dotsb
\xrightarrow{f_{n+1}} X_{n+1}
\]
is a composite of $n+1$ ghosts.  Each of the maps $f_{i}\mathcolon
X_{i}\xrightarrow{}X_{i+1}$ has positive filtration.  Therefore, their
composition will have filtration $\geq n+1$ by Proposition~IV.4.4
of~\cite{elmendorf-kriz-mandell-may}, and therefore will be null.
Conversely, suppose $X$ has $\pdim X\leq n$.  In the construction of
the universal coefficient spectral sequence, we construct exact
triangles 
\begin{gather*}
\Sigma^{-1} X_{1} \xrightarrow{} P_{0} \xrightarrow{g_{0}} X
\xrightarrow{h_{1}} X_{1} \\
\Sigma^{-1}X_{2} \xrightarrow{} P_{1} \xrightarrow{g_{1}} X_{1}
\xrightarrow{h_{2}} X_{2} \\
\dotsb \\
\Sigma^{-1} X_{j+1} \xrightarrow{} P_{j} \xrightarrow{g_{j}} X_{j}
\xrightarrow{h_{j+1}} X_{j+1}
\end{gather*}
where $P_{j}$ is projective and $\pi_{*}g_{j}$ is onto for all $j$,
and so $h_{j+1}$ is a ghost for all $j$.  A map has filtration $n+1$
if and only if it factors through the composite
\[
X \xrightarrow{h_{1}} X_{1} \xrightarrow{h_{2}} \dotsb
\xrightarrow{h_{n+1}} X_{n+1}.
\]
Thus every map of filtration $n+1$ is a composite of $n+1$ ghosts.
Thus, $\pdim X\leq n$ implies every map of filtration $n+1$ is null.  

Part~(2) is then completely dual.  Part~(3) is proved
in~\cite{hovey-lockridge-weak}.  
\end{proof}

The following fact is also useful. 

\begin{proposition}\label{prop-alg-dim}
Let $E$ be an $S$-algebra, and let $X$ be an object of $\cat{D} (E)$.
Then
\[
\pdim X\leq \pdim_{E_{*}} X_{*}, \ \ \idim X\leq \idim_{E_{*}} X_{*},
\]
and 
\[
\fdim X\leq \cfdim X \leq \fdim_{E_{*}} X_{*}.  
\]
\end{proposition}

The first two statements of this proposition are obvious given
Proposition~\ref{prop-proj-dim}, and the last statement is proved
in~\cite{hovey-lockridge-weak}.

We can now define our homological dimensions. 

\begin{definition}\label{defn-global}
Suppose $E$ is an $S$-algebra.  Define the \textbf{(right) global
dimension} of $E$, $\rgldim E$, by 
\[
\rgldim E = \sup_{X} \pdim X =\sup_{X} \idim X.
\]
These two numbers are equal because they are both equal to the longest
nontrivial composition of ghosts in $\cat{D} (E)$ (or $\infty$ if
there are arbitrarily long nontrivial compositions of ghosts).  In
case $E$ is a commutative $S$-algebra, we just refer to the
\textbf{global dimension}, $\gldim E$.  Define
the \textbf{ghost dimension} of $E$, $\ghdim E$, by
\[
\ghdim E = \sup_{X \text{compact}} \pdim X =\sup_{X} \fdim X. 
\]
\end{definition}

It is proven in~\cite{hovey-lockridge-weak} that the two definitions
of ghost dimension given above coincide.  

We then have the following theorem that sums up the basic properties
of these definitions.  

\begin{theorem}\label{thm-ring}
Suppose $E$ is an $S$-algebra.
\begin{enumerate}
\item $\ghdim E\leq \rgldim E$.  
\item $\rgldim E\leq \rgldim E_{*}$, with equality if $E=HR$ for an
ordinary ring $R$.  
\item $\ghdim E\leq \wdim E_{*}$, with equality if $E=HR$ for an
ordinary ring $R$.
\item $\ghdim E=\ghdim E^{\textup{op}}$, where $E^{\textup{op}}$ is
$E$ with the opposite multiplication.  
\end{enumerate}
\end{theorem}

The first part of this theorem is obvious.  Part~(2) is proved
in~\cite{hovey-lockridge-ghost}, though it is originally in the second
author's thesis.  Parts~(3) and~(4) are proved
in~\cite{hovey-lockridge-weak}.  

$S$-algebras of global dimension $0$ are called \textbf{semisimple}.
They are studied in~\cite{hovey-lockridge-semisimple}.  In particular,
there are semisimple ring spectra with $\rgldim E_{*}=\infty$.
Similarly, $S$-algebras of ghost dimension $0$ are called \textbf{von
Neumann regular}, and are also studied
in~\cite{hovey-lockridge-semisimple}.

\section{Examples}\label{sec-examples}

We would, of course, like to compute $\ghdim E$ and $\rgldim E$ for
various $S$-algebras $E$.  Some of this was done
in~\cite{hovey-lockridge-semisimple}, where the authors classified the
semisimple $S$-algebras $E$ if either $E_{*}$ is commutative or local.
Those spectra are rather unusual, however.  We address some of the
more common $S$-algebras $E$ in this section, concentrating on the
case when $E_{*}$ is Noetherian.

We begin with the sphere spectrum, and the following unsurprising
result.  

\begin{proposition}\label{prop-sphere}
Let $S$ be the sphere $S$-algebra.  Then $\ghdim S=\gldim S=\infty$.  
\end{proposition}

\begin{proof}
This is due to Christensen~\cite{christensen}, who provides bounds on
$\pdim \mathbb{R}P^{n}$ (which he calls the length).  In particular, a
lower bound for $\pdim \mathbb{R}P^{k}$ is given by the longest
nonzero chain of Steenrod operations in its homology, since Steenrod
operations are obviously ghosts.  And this longest chain is easily
seen to grow without bound as $k$ grows.
\end{proof}

\begin{corollary}\label{cor-sphere-local}
Let $S_{(p)}$ be the $p$-local sphere $S$-algebra, where $p$ is an
integer prime.  Then $\ghdim S_{(p)}=\gldim S_{(p)}=\infty$.  
\end{corollary}

\begin{proof}
Use reduced power operations in the cohomology of skeleta of $B\Z /p$
to replace Steenrod operations in $\mathbb{R}P^{k}$.  We do not need
to know the exact length of a nontrivial composition of these
operations, we just need to know that this length grows without bound
as we take larger skeleta. 
\end{proof}

Recall that the length of the longest regular sequence in a ring is
often called the \textbf{depth}.   

\begin{theorem}\label{thm-regular}
If $E$ is a commutative $S$-algebra, then 
\[
\depth E_{*} \leq \ghdim E \leq \min \{\wdim E_{*}, \rgldim E \} \leq
\rgldim E_{*}.
\]
\end{theorem}

\begin{proof}
In view of Theorem~\ref{thm-ring}, it suffices to prove the first
inequality.  Let 
\[
x_{1},x_{2},\dotsc ,x_{n}
\]
be a regular sequence in $E_{*}$.  The main algebraic input we need is
the computation that
\[
\Ext^{i}_{E_{*}} (E_{*}/ (x_{1},\dotsc ,x_{n}),E_{*})=0 \text{ if } i\neq n
\]
and is nonzero if $i=n$.  One can prove this by induction on $n$ (or
by the Koszul resolution), using the exact sequences
\[
0 \xrightarrow{} E_{*}/ (x_{1},\dotsc ,x_{k-1}) \xrightarrow{x_{k}}
E_{*}/ (x_{1},\dotsc ,x_{k-1}) \xrightarrow{} E_{*}/ (x_{1},\dotsc
,x_{k}) \xrightarrow{} 0,
\]
where we have ignored suspensions for simplicity. 

We also need the fact that there is a $E$-module $E/ (x_{1},\dotsc
,x_{n})$ realizing the $E_{*}$-module $E_{*}/ (x_{1},\dotsc ,x_{n})$.
One can also construct these by induction, using the exact triangles
\[
E/ (x_{1},\dotsc ,x_{i-1})\xrightarrow{x_{i}} E/ (x_{1},\dotsc
,x_{i-1}) \xrightarrow{} E/ (x_{1},\dotsc ,x_{i}) \xrightarrow{}  E/
(x_{1},\dotsc ,x_{i-1}),
\]
ignoring suspensions again.  

The universal coefficient spectral sequence 
\[
\Ext_{E_{*}}^{s,t} (E_{*}/ (x_{1},\dotsc ,x_{n}),E_{*})\Rightarrow
\cat{D} (E) (E/ (x_{1},\dotsc ,x_{n}),E) 
\]
then has only one non-vanishing line, where $s=n$.  It therefore
collapses, and so there is an element in $E_{\infty}$ of filtration
$n$.  Thus $\ghdim E\geq n$, as required.  
\end{proof}

It is tempting to believe that Theorem~\ref{thm-regular} works in the
noncommutative case as well, as long as we take regular sequences in
the center $Z (E_{*})$.  The algebraic calculation works fine, but we
do not seem to be able to construct the necessary maps $x\mathcolon
M\xrightarrow{}M$ for an $E$-module $M$ and an $x\in Z (E_{*})$.  We
can construct such maps for $E$-bimodule maps $x\mathcolon
E\xrightarrow{}E$, but an element in the center of $E_{*}$ need not
give such a map, so far as we know.  

\begin{corollary}\label{cor-MU}
We have 
\[
\ghdim MU =\ghdim BP =\infty,
\]
while $\ghdim E_{n}=\gldim E_{n}=n$ and $\ghdim K=\gldim K=1$.
\end{corollary}

Here $E_{n}$ denotes Morava $E$-theory, with 
\[
E_{n*}\cong W\mathbb{F}_{p^{n}}[[u_{1},\dotsc ,u_{n-1}]][u,u^{-1}].  
\]
This is known to be a commutative $S$-algebra
by~\cite{goerss-hopkins-rings}.  The coefficient ring has global
dimension $n$, and $p,u_{1},\dotsc ,u_{n-1}$ is a regular sequence,
so Theorem~\ref{thm-regular} gives the desired result.  

This corollary suggests the following more general theorem.

\begin{theorem}\label{thm-comm-Noeth}
Suppose $E$ is a commutative $S$-algebra such that $E_{*}$ is
Noetherian with $\gldim \pi_{*}E<\infty$.  Then $\ghdim E=\gldim
E=\gldim E_{*}$.
\end{theorem}

\begin{proof}
The point is that $\depth R=\gldim R$ when $R$ is Noetherian
commutative of finite global dimension, so the result follows from
Theorem~\ref{thm-regular}.  This algebraic fact is a corollary of
Serre's characterization of regular local rings, but we
cite~\cite[Section~2]{brown-hajarnavis} because it contains an
interesting non-commutative generalization as well.  
\end{proof}

We suspect that Theorem~\ref{thm-comm-Noeth} may be true even if $E$
is not a commutative $S$-algebra.  Much less is known about
noncommutative Noetherian rings of finite global dimension, however.


We now address the global dimension of real $K$-theory.
Unfortunately, we have not been able to determine the exact value of
$\gldim KO$, but we do at least bound it.  

\begin{theorem}\label{thm-KO}
Let $KO$ denote $2$-local periodic real $K$-theory and $ko$ denote
$2$-local connective real $K$-theory, both of which are commutative
$S$-algebras.  Then
\[
1\leq \gldim KO \leq 3
\]
and 
\[
4\leq \gldim ko \leq 5.
\]
\end{theorem}

\begin{proof}
This depends on the results of Bousfield~\cite{bousfield-united} and
Wolbert~\cite{wolbert}.  Bousfield wrote his paper using naive
$KO$-module spectra, as opposed to objects in $\cat{D} (KO)$.
However, Wolbert explains why Bousfield's results hold in $\cat{D}
(KO)$ as well, and proves analogous results for $\cat{D} (ko)$.  We
note that there is an error in Wolbert's paper~\cite{sagave-toda}, but
this error is about the realizability of $ko_{*}$ and $KO_{*}$-modules
and does not affect this theorem.

Bousfield constructs an abelian category of CRT-modules, which are
modules over the 3-object additive category consisting of
$\{KO_{*},K_{*}, KSC_{*} \}$ and the various standard maps between
them.  Here $KSC=KO\wedge C (\eta^{2})$ is self-conjugate $K$-theory.
There is a functor from $KO$-module spectra to this category that
takes $X$ to $\pi_{*}^{CRT} (X)$, which is the set $\{\pi_{*}X,
\pi_{*} (C (\eta)\wedge X), \pi_{*} (C (\eta^{2})\wedge X) \}$
together with the maps between them.  He then proves that
$\pi_{*}^{CRT} (X)$ has projective dimension $\leq 1$ for every
$KO$-module spectrum $X$.  The $KO$-modules $P$ such that
$\pi_{*}^{CRT}$ is projective as a CRT-module are coproducts of
suspensions of $KO, K$, and $KSC$, and every projective CRT-module
arises this way.  Furthermore, maps between projective CRT-modules are
realizable as maps of the corresponding $KO$-modules.  This means that
for every $KO$-module spectrum $X$, there is a cofiber sequence
\[
Q \xrightarrow{} P \xrightarrow{} X \xrightarrow{} \Sigma Q
\]
in $\cat{D} (KO)$ in which $Q$ and $P$ are coproducts of suspensions
of copies of $KO,K$, and $KSC$.  Since $K=KO\wedge C (\eta)$ and
$KSC=KO\wedge C (\eta^{2})$, we see that any $2$-fold composite of
ghosts out of $P$ or $Q$ is trivial.  A simple argument then shows
that any $4$-fold composite of ghosts out of $X$ is trivial, so
$\gldim KO\leq 3$.  

The argument for $ko$ is a little more
complicated. Wolbert~\cite{wolbert} describes the connective version
of crt-modules, and shows that if $X$ is a $ko$-module, then
$\pi_{*}^{crt} (X)$ has projective dimension at most $2$ in the
category of crt-modules.  Wolbert then gives a version of the
universal coefficient spectral sequence 
\[
E^{s,t}_{2} = \Ext_{crt}^{s,t} (\pi_{*}^{crt} (X), \pi_{*}^{crt} (Y))
\Rightarrow \cat{D} (ko) (X,Y)_{t-s}. 
\]
As before, any $2$-fold composite of ghosts has filtration $1$ in this
spectral sequence (since $k$ and $ksc$ are $2$-cell complexes in
$\cat{D} (ko)$).  Therefore, any $6$-fold composite of ghosts will
have filtration $3$, but the spectral sequence is trivial above
filtration $2$.  Thus every $6$-fold composite of ghosts is trivial,
so $\gldim ko\leq 5$.  

To see that $\ghdim ko\geq 4$, we use the fact that $ko\wedge A
(1)=H\mathbb{F}_{2}$, where $A (1)$ is the usual $8$-cell complex with
whose cohomology is $\mathcal{A}<x>/ (\Sq^{2^{n}}x|n\geq 2)$.  Using
this, we can compute $\cat{D} (ko)
(H\mathbb{F}_{2},H\mathbb{F}_{2})$.  It is the subring of the Steenrod
algebra generated by $\Sq^{1}$ and $\Sq^{2}$.  Since every element of
the Steenrod algebra is a ghost, the nontrivial element
$\Sq^{2}\Sq^{1}\Sq^{2}\Sq^{1}$ is a nontrivial composite of $4$
ghosts, and is a self-map of the compact object $H\mathbb{F}_{2}$ in
$\cat{D} (ko)$.  Thus $\ghdim ko\geq 4$.  
\end{proof}

The proof of this theorem is of course dependent on knowing an awful
lot about $KO$ and $ko$.  We have much less information about higher
analogues of $KO$, such as the spectrum $tmf$ of topological modular
forms~\cite{hopkins-tmf}.  However, it is often the case with such
spectra $E$ that there is a finite type $0$ spectrum $X$ such that
$E\wedge X$ is a Noetherian $S$-algebra of finite global dimension.
For example, $KO\wedge C (\eta)=KU$, and at the prime $3$, there is a
3-cell complex $T$ such that $tmf\wedge T$ is a wedge of two copies of
$BP<2>$~\cite[Lemma~2, after Corollary~2.4.6]{behrens-k2}.  Presumably
a larger such finite complex also exists at the prime $2$ for $tmf$,
though such a result has not been proven as yet.  

In general, given spectra $X$ and $Y$, we can define the term ``$Y$
can be built from $X$ in $\ell$ steps'' in the same way that we
defined the projective dimension.  That is, we
say that $Y$ can be built from $X$ in $0$ steps if $Y$ is a retract of
a coproduct of suspensions of $X$.  We then say that
$Y$ can be built from $X$ in $\ell$ steps if there is an exact triangle
\[
Z \xrightarrow{} W \xrightarrow{} \widetilde{Y} \xrightarrow{} \Sigma Z
\]
where $W$ can be built from $X$ in $0$ steps, $Z$ can be built from
$X$ in $\ell -1$ steps, and $Y$ is a retract of $\widetilde{Y}$.  

\begin{theorem}\label{thm-finite}
Suppose $E$ is an $S$-algebra and $X$ is a spectrum with the following
properties.  
\begin{enumerate}
\item $E\wedge X$ is an $E$-algebra, so also an $S$-algebra, with
$\rgldim (E\wedge X)=m<\infty$.
\item As an object of $\cat{D} (S)$, $\pdim X=k$.  
\item $S$ can be built from $X$ in $\ell$ steps.  
\end{enumerate}
Then $\gldim E\leq (k+1) (\ell +1) (m+1)-1$.  
\end{theorem}

Note that if $Y$ is in the thick subcategory generated by $X$, then
$Y$ can be built from $X$ in a finite number of steps.  In particular,
if $X$ is a type $0$ finite spectrum, then $S$ can be built from $X$
in a finite number of steps.  Hence we get the following corollary.  

\begin{corollary}\label{cor-finite}
Suppose $E$ is an $S$-algebra and $X$ is a type $0$ finite spectrum
such that $E\wedge X$ is an $E$-algebra with finite right global
dimension as an $S$-algebra.  Then $E$ has finite right global
dimension.
\end{corollary}

Here is the proof of Theorem~\ref{thm-finite}. 

\begin{proof}
Note that if $M$ is an $E$-module, then $M\wedge X\cong M\wedge_{E}
(E\wedge X)$ is an $E\wedge X$-module.  Similarly, if $f$ is a map of
$E$-modules, then $f\wedge 1_{X}$ is a map of $E\wedge X$-modules.  

Now, let $M$ and $N$ be $E$-modules.  By induction on $t$, one can see
that if $\pdim Z=t$, then every $t+1$-fold ghost $f\mathcolon
M\xrightarrow{}N$ induces a ghost $f\wedge 1_{Z}\mathcolon M\wedge
Z\xrightarrow{}N\wedge Z$ of $E$-modules.  Taking $Z=X$, we see that
if $f$ is a $k+1$-fold ghost, then $f\wedge 1_{X}$ is a ghost,
necessarily as a map of $E\wedge X$-modules. Hence, if $f\mathcolon
M\xrightarrow{}N$ is a $(k+1) (m+1)$-fold ghost, then $f\wedge 1_{X}$
is an $(m+1)$-fold ghost, and hence is null as a map of $E\wedge
X$-modules, and in particular as a map of $E$-modules.  But then we
can proceed by induction on $\ell$ to see that if $Y$ can be built
from $X$ in $\ell$ steps, then any $(k+1) (m+1) (\ell +1)$-fold ghost
$f$ has $f\wedge 1_{Y}$ null as a map of $E$-modules.  Taking $Y=S$
completes the proof.
\end{proof}

\begin{corollary}\label{cor-tmf}
At the prime $3$, the spectrum $tmf$ of topological modular forms has
finite global dimension.
\end{corollary}

\begin{proof}
As mentioned above, there is a $3$-cell complex $T$ such that
\[
tmf\wedge T\simeq BP<2> \wedge \Sigma^{8} BP<2>
\]
by~\cite{behrens-k2}.  The complex $T$ has cells in dimensions $0$,
$4$, and $8$, so is obviously type $0$.  The spectrum $tmf\wedge T$ is
also called $tmf_{0} (2)$, and is a commutative $tmf$-algebra (it is
the connective cover of $TMF_{0} (2)$, which is the spectrum of
sections over a certain stack of a sheaf of commutative $S$-algebras;
see~\cite{behrens-k2}).  The homotopy ring of $tmf_{0} (2)$ is
polynomial over $\Z_{3}$ (if we use the completed version) on two
generators~\cite[Proposition~2.3]{hill-tmf}, and therefore $tmf_{0}
(2)$ has global dimension $3$ by~\ref{thm-comm-Noeth}.
\end{proof}

\section{Properties}\label{sec-properties}

In this section, we examine some general properties of global and
ghost dimension.  Most of these properties concern the relationship
between the global dimension or ghost dimension of an $S$-algebra $E$
and some other $S$-algebra $F$ related to it.  We discuss the cases
when $F$ is a smashing Bousfield localization of $E$, when $F$ is a
free $E$-module, and when $F$ is Morita equivalent to $E$.  We end
with a discussion of the relationship between ghost and global
dimension.  We find this particularly unsatisfactory, however, because
we are unable to prove anything along the lines of the well-known
algebraic fact that if $R$ is Noetherian or right perfect (which is
equivalent to flats being projective), then $\rgldim R=\wdim R$.  

\begin{proposition}\label{prop-localization}
Suppose $L$ is a smashing Bousfield localization functor, and $E$
is an $S$-algebra.  Then 
\[
\rgldim E\geq \rgldim LE.
\]
\end{proposition}

We do not know if this theorem is true if the localization functor is
not smashing.  

\begin{proof}
Note first that $LE$ is again an
$S$-algebra~\cite[Chapter~VIII]{elmendorf-kriz-mandell-may}.  The main
point is that because $L$ is smashing, the category $L\cat{D}
(E)$ of $L$-local $E$-modules is equivalent to the category
$\cat{D}
(LE)$~\cite[Proposition~VIII.3.2]{elmendorf-kriz-mandell-may}.
Thus a nontrivial composite of ghosts in $\cat{D} (LE)$ is the
same thing as a nontrivial composite of ghosts between $L$-local objects
in $\cat{D} (E)$.  Hence $\rgldim E\geq \rgldim LE$.  
\end{proof}

The same thing is true for ghost dimension, though we need to assume
$\cat{D} (E)$ is a Brown category.  Recall from Section~4.2
of~\cite{hovey-axiomatic} that this means homology theories, and
morphisms between them, are representable.  

\begin{proposition}\label{prop-localization-ghost}
Suppose $L$ is a smashing Bousfield localization functor, and $E$ is
an $S$-algebra such that $\cat{D} (E)$ is a Brown category.  Then 
\[
\ghdim E \geq \ghdim LE.
\]
\end{proposition}

\begin{proof}
Because $\cat{D} (E)$ is a Brown category, every object of $\cat{D}
(E)$ is the minimal weak colimit of the compact objects mapping to
it by~\cite[Theorem~4.2.4]{hovey-axiomatic}.  But then we can follow
the proof of~\cite[Theorem~6.2 (b,c)]{hovey-strickland} to show that
the compact objects of $L\cat{D} (E)$ are the retracts of objects of
the form $LF$, for $F$ compact in $\cat{D} (E)$.  Thus, if we have a
nontrivial composite of ghosts 
\[
X_{0} \xrightarrow{} X_{1} \xrightarrow{} \dotsb \xrightarrow{} X_{n}
\]
in $L\cat{D} (E)$, where $X_{0}$ is compact in $L\cat{D} (E)$, we can
write $X_{0}$ as a retract of $LF$ for some compact object $F\in
\cat{D} (E)$.  Then the composite 
\[
F \xrightarrow{} LF \xrightarrow{} X_{0}\xrightarrow{} X_{1}
\]
must be nontrivial, and gives us a nontrivial composite of ghosts out
of a compact object in $\cat{D} (E)$.  
\end{proof}

In general, there is not much relationship between the global
dimension of an $S$-algebra $E$ and a general $E$-algebra $F$.  At one
extreme we have the smashing localizations discussed above.  The other
extreme is $E$-algebras which are free over $E$.  

\begin{proposition}\label{prop-free}
Suppose $E\xrightarrow{}F$ is a map of $S$-algebras such that $F_{*}$
is free over $E_{*}$.  Then $\rgldim E\leq \rgldim F$ and $\ghdim E\leq
\ghdim F$.  
\end{proposition}

Note that the inequalities in this proposition can certainly be
strict, as we can see by looking at the inclusion of ordinary rings
from $\Z$ to $\Z [x]$, for example.  We think this proposition should
hold even if $F_{*}$ is only assumed to be faithfully flat over
$E_{*}$, but have not been able to prove it.

\begin{proof}
Consider the extension functor $F\wedge_{E} (-)\mathcolon \cat{D}
(E)\xrightarrow{}\cat{D} (F)$ and its right adjoint, the restriction
functor.  Because $F_{*}$ is flat over $E_{*}$, the natural map
\[
F_{*}\otimes_{E_{*}} X_{*} \xrightarrow{} \pi_{*} (F\wedge_{E}X)
\]
is an isomorphism.  Indeed, both sides are homology functors on
$\cat{D} (E)$, and the given map is an isomorphism for $X=E$, so it is
an isomorphism for all $X\in \cat{D} (E)$.  This means that
$F\wedge_{E} (-)$ preserves ghosts.  

In particular, if $g$ is a composition of $n$ ghosts, and $F_{*}$ is
flat over $E_{*}$, then $F\wedge_{E}g$ is a composition of $n$
ghosts.  If the domain of $g$ is a compact object of $\cat{D} (E)$,
then the domain of $F\wedge_{E}g$ is a compact object of $\cat{D}
(F)$.  However, $F\wedge_{E}g$ may be zero, even if $g$ is nonzero.  This
cannot happen if $F_{*}$ is free over $E_{*}$, however, for then
$F\wedge_{E}X$ is a coproduct of copies of $X$ as an $E$-module.  Thus
$g$ is a retract of the restriction of $F\wedge_{E}g$, so if $g$ is
nontrivial, so is $F\wedge_{E}g$.   
\end{proof}

We now discuss Morita equivalence. Recall from the work of Scwede and
Shipley~\cite{schwede-morita} that two $S$-algebras $E$ and $F$ are
called Morita equivalent if there is a chain of Quillen equivalences
from the model category of $E$-modules to the model category of
$F$-modules.  Schwede and Shipley prove that, if $E$ and $F$ are
Morita equivalent and cofibrant in the model structure on $S$-algebras
(we can always assume this, since weak equivalences of $S$-algebras
induce Quillen equivalences of the module categories), then there is
an $E-F$-bimodule $M$ and an $F-E$-bimodule $N$, both of which are
compact both as $E$ and $F$-modules, so that the functors
\[
\Phi (X) = X\wedge_{E} M \text{ and } \Psi (Y) = Y\wedge_{F} N
\]
are inverse equivalences, with $\Phi \mathcolon \cat{D}
(E)\xrightarrow{}\cat{D} (F)$ and $\Psi$ going back the other way.  
Furthermore, $M$ generates $\cat{D} (F)$ and $N$ generates $\cat{D}
(E)$, in the sense that the smallest localizing subcategory containing
$M$ (resp. $N$) is $\cat{D} (E)$ (resp. $\cat{D} (F)$).  

This is analogous to the usual Morita situation with ordinary rings,
with one important difference.  The generators $M$ and $N$ do not have
to be projective.  This means that neither $\Phi$ nor $\Psi$ need
preserve projective objects, so that we do NOT expect global dimension
or ghost dimension to be Morita invariant.

Indeed, we thank Lidia Angeleri H\"{u}gel for the following example
coming from the theory of tilting modules, a reference for which
is~\cite{assem-simson-skowronski}.  Recall that, for an ordinary ring
$A$, a tilting module is a module $T$ whose projective resolution
gives an equivalence of categories from the derived category of $A$ to
the derived category of $\End_{A} (T)$ by using the derived tensor
product.  In particular, a tilting module defines a Morita equivalence
in the above sense between the Eilenberg-MacLane $S$-algebras $HA$ and
$HB$.  Now, let $A$ be the path algebra of the quiver $1
\xrightarrow{}2 \xrightarrow{}3$, so that $A$ is isomorphic to the
ring of lower triangular $3\times 3$ matrices.  There are then three
projective indecomposable $A$-modules $P_{1},P_{2}$, and $P_{3}$,
corresponding to the quiver representations
\[
k \xrightarrow{=}k \xrightarrow{=}k, \ \
0\xrightarrow{}k\xrightarrow{=}k, \text{ and }
0\xrightarrow{}0\xrightarrow{}k.
\]
There are also three dual injective indecomposables
$I_{1},I_{2},I_{3}=P_{1}$ corresponding to the representations
\[
k\xrightarrow{}0 \xrightarrow{}0, \ \ k\xrightarrow{=}k\xrightarrow{}0,
\text{ and } k\xrightarrow{=}k\xrightarrow{=}k.
\]

We claim that $T=P_{3}\oplus I_{1}\oplus P_{1}$ is a tilting module.
Indeed, like any path algebra, $A$ has global dimension $\leq 1$ (and
in fact $\rgldim A=1$), so $\pdim T\leq 1$.  We also need
$\Ext^{1}(T,T)=0$, and this boils down to $\Ext^{1} (I_{1},P_{3})=0$,
which is straightforward.  The last condition for $T$ to be a tilting
module is for there to be an exact sequence
\[
0 \xrightarrow{} A \xrightarrow{} T_{1} \xrightarrow{} T_{2}
\xrightarrow{} 0
\]
where $T_{1}$ and $T_{2}$ are finite direct sums of direct summands of
$T$, so finite direct sums of $P_{3}$, $I_{1}$, and $P_{1}$.  But $A$
corresponds to the representation 
\[
k \xrightarrow{} k\oplus k \xrightarrow{} k\oplus k\oplus k
\]
where the maps are inclusions of the obvious summands.  Thus we can
just take $T_{2}=I_{1}$ and $T_{1}=P_{3}\oplus P_{1}\oplus P_{1}$ with
the obvious surjection from $T_{1}$ to $T_{2}$.  

So we have a Morita equivalence between $HA$ and $HB$, where
\[
B=\End_{A} (T).
\]
One can compute $B$ directly.  It is a
$5$-dimensional algebra generated by the orthogonal idempotents
$e_{1}, e_{2}, e_{3}$ and two other elements $\alpha ,\beta$, where
the only nonzero products involving $\alpha$ and $\beta$ are
\[
e_{3}\alpha =\alpha e_{1}=\alpha  \text{ and } e_{2}\beta =\beta e_{3}=\beta .
\]
Let $M$ be the right $B$-module which has dimension $1$ over $k$ where
$e_{2}$ acts by the identity and the other generators of $B$ act
trivially.  Then one can check that $\pdim M=2$, so 
\[
\rgldim B\geq 2>\rgldim A = 1.  
\]
In fact, $\rgldim B=2$.  Note that both $A$ and $B$ are Noetherian,
so their weak dimensions are also different.  Thus $HA$ and $HB$ have
different ghost dimensions as well.  

We should mention that the property of having global dimension $0$ IS
Morita invariant, at least if one of the two $S$-algebras has
commutative
homotopy~\cite[Proposition~2.11]{hovey-lockridge-semisimple}.  We now
prove that the finiteness of global dimension is at least Morita
invariant.  

\begin{proposition}\label{prop-Morita}
Suppose $E$ and $F$ are Morita equivalent $S$-algebras.  Then $\rgldim
E<\infty$ if and only if $\rgldim F<\infty$.  Similarly, $\ghdim
E<\infty$ if and only if $\ghdim F<\infty$.
\end{proposition}

For example, this means that the endomorphism $S$-algebra of any
finite type $0$ spectrum has infinite global dimension.  In general,
we expect $\rgldim E$ to be always infinite if $E$ is a finite
spectrum.  

\begin{proof}
As mentioned above, we can assume that $E$ and $F$ are cofibrant
$S$-algebras, and we have compact generators $M$ of $\cat{D} (F)$ and
$N$ of $\cat{D} (E)$ that are also bimodules, so that smashing with
$M$ and $N$ give inverse equivalences.  Now suppose $f$ is a ghost map
in $\cat{D} (E)$. Then $f\wedge_{E}M$ has the property that $\cat{D}
(F) (M,f\wedge_{E}M)_{*}=0$.  Since $M$ is a compact generator of
$\cat{D} (F)$, $F$ is in the thick subcategory generated by $M$.
Thus, $F$ can be built from $M$ in finitely many steps (see the
discussion before Theorem~\ref{thm-finite}).  In particular, there
is an integer $k$ such that, if $f$ is a composite of $k$ ghosts in
$\cat{D} (E)$, then $f\wedge_{E}M$ is a ghost in $\cat{D} (F)$.  In
particular, if $F$ has finite global dimension, say $n$, then if $f$
is a composite of $k (n+1)$ ghosts in $\cat{D} (E)$, then
$f\wedge_{E}M$ is null.  Since smashing with $M$ is an equivalence of
categories, this means that $f$ is null, so $E$ has finite global
dimension.  

Reversing the roles of $E$ and $F$ gives the reverse implication.  For
the ghost dimension, we repeat the same argument using a compact
object as the source of our first ghost map.  This causes no problems
since the equivalences of categories preserve compact objects.  
\end{proof}

We now look for a relationship between the global dimension and the
ghost dimension.  In algebra, we have the well-known inequality 
\[
\rgldim R \leq \wdim R + \pgldim R,
\]
where $\pgldim R$ is the pure global dimension of $R$.  This can
actually be replaced by the maximum projective dimension of a flat
module.  Indeed, in a projective resolution of an arbitrary $R$-module
$M$, the syzygies $M_{k}$ are flat whenever $k\geq \wdim R$.  But this
means our projective resolution is also a projective resolution of the
flat module $M_{\wdim R}$, giving us the desired inequality.

One might expect the analogue of the pure global dimension to be
the \emph{phantom dimension}, defined below.  

\begin{definition}\label{defn-phantom}
Suppose $E$ is an $S$-algebra.  Define the \textbf{phantom dimension}
of $E$, $\phdim E$, to be the smallest integer $n$ such that every
composite of $n+1$ phantom maps in $\cat{D} (E)$ is zero, or $\infty$
if there exist arbitrarily long such nonzero composites.
\end{definition}

Recall that a phantom map is a map $f$ for which $[C,f]_{*}=0$ for
every compact $E$-module $C$.  So, if we think of a map whose cofiber
is a ghost as the homotopy-theoretic analogue of an epimorphism, then a
map whose cofiber is a phantom is the homotopy-theoretic analogue of a
pure epimorphism, and the phantom dimension should be analogous to the
pure global dimension.  Furthermore, the phantom dimension is
obviously invariant under Morita equivalences, since it only mentions
compact objects, which are preserved by any equivalence of
categories.  

Note that if $E_{*}$ is countable, or, more generally, if $\cat{D}
(E)$ is a Brown category, then $\phdim E\leq
1$~\cite{neeman-brown,christensen}.  Very little else about the
phantom dimension is known, except that it can be greater than one (as
is shown in the above papers).  The natural conjecture is that $\phdim
HR$ should be the pure global dimension of $R$, but this is false for
$R=\Z /4$, whose pure global dimension is $0$ whereas $\phdim \Z
/4=1$.  The problem here is the difference between finitely presented
$R$-modules and compact objects of $\cat{D} (R)$.

In any case, we have the following proposition. 

\begin{proposition}\label{prop-phantom}
Suppose $E$ is an $S$-algebra.  Then 
\[
\rgldim E <(\ghdim E +1) (\phdim E+1). 
\]
In particular, if $\cat{D} (E)$ is a Brown category, then 
\[
\rgldim E \leq 2\ghdim E + 1.  
\]
\end{proposition}

This proposition is quite a bit weaker than the inequality for
ordinary rings, mentioned above.  This is because we cannot construct
a resolution of an $E$-module in the same way as we can in algebra.
But possibly this proposition can be improved.  

\begin{proof}
Suppose $\ghdim E=n$ and $\phdim E=k$.  There is nothing to prove if
either of these is infinite, so assume they are finite.  Then any
$n+1$-fold composite of ghosts is phantom, so any $(n+1) (k+1)$-fold
composite of ghosts is a composite of $k+1$ phantom maps, and so is
null.  
\end{proof}



\section{Gorenstein rings}\label{sec-Gorenstein}

Recall that one of the themes of~\cite{hovey-lockridge-semisimple} was
that if $E$ is an $S$-algebra and $\rgldim E=0$, then there are very
severe restrictions on $E_{*}$.  In particular, we show that $E_{*}$
must be quasi-Frobenius, though much more is true.  Since
quasi-Frobenius rings are the same as $0$-Gorenstein rings, a natural
conjecture might be that if $\rgldim E=n$, then $E_{*}$ is
$n$-Gorenstein.  Unfortunately, this is easily seen to be false, since
$KO_{*}$ is not $n$-Gorenstein for any $n$, as we show in this
section.  However, we also show that if $\rgldim E=1$ and $E_{*}$ is a
Noetherian domain, then $E_{*}$ is $1$-Gorenstein.  We do not know if
this statement is true for larger $n$, though it seems unlikely.  

Recall that a (possibly noncommutative) ring $R$ is called
\textbf{Gorenstein} if it is left and right Noetherian and $R$ has
finite injective dimension as a left or right $R$-module.  This
generalizes the usual definition of Gorenstein in the commutative
case, which is much used in algebraic geometry.  If $R$ is Gorenstein,
the right and left injective dimensions of $R$ must coincide, and if
they are at most $n$, $R$ is called $n$-Gorenstein.  These
generalizations of quasi-Frobenius rings, which are just
$0$-Gorenstein rings, have been the object of much recent
study. Chapter~9 of~\cite{enochs-jenda-book} is a good place to start.

We begin by showing that $KO_{*}$ is not Gorenstein. 

\begin{proposition}\label{prop-KO}
The ring $KO_{*}$ is not $n$-Gorenstein for any $n$, though it is
Noetherian and $\gldim KO<\infty$.  
\end{proposition}

\begin{proof}
If $KO_{*}$ were Gorenstein, it would be Cohen-Macaulay, which would mean
that its depth would be equal to its Krull dimension.  But 
\[
KO_{*} = \Z_{(2)}[\eta ,w,v,v^{-1}]/ (\eta^{3}, 2\eta , w\eta , w^{2}-4v).
\]
The only prime ideas in $KO_{*}$ are $(\eta)$ and the maximal ideal
$(\eta ,2,w)$.  So the Krull dimension is $1$.  But there are no
non-zero divisors in the maximal ideal, so the depth is $0$, and so
$KO_{*}$ is not Cohen-Macaulay.  
\end{proof}

We now consider the case when $\gldim E=1$.  

\begin{proposition}\label{prop-dim-one}
Suppose $E$ is an $S$-algebra with 
\[
\rgldim E=\rgldim E^{\textup{op}}=1,
\]
and suppose that there are no nonzero maps from an
injective \ulp left or right\urp \ $E_{*}$-module to a projective \ulp
left or right\urp \ $E_{*}$-module.  Then $\idim E_{*}\leq 1$ as
either a left or right $E_{*}$-module. Consequently, if $E_{*}$ is
also left and right Noetherian, then $E_{*}$ is a $1$-Gorenstein ring.
\end{proposition}

\begin{proof}
Let $R=E_{*}$.  It suffices to show that $\idim R\leq 1$.  We just do
this on one side, since the hypotheses are left-right symmetric.  For
this, we embed $R$ into an injective module $I_{0}$, giving us the
short exact sequence
\[
0 \xrightarrow{} R \xrightarrow{} I_{0} \xrightarrow{} R_{1}
\xrightarrow{} 0. 
\]
We can realize this uniquely by an exact triangle in $\cat{D} (E)$
\[
E \xrightarrow{} J_{0} \xrightarrow{} M_{1} \xrightarrow{\delta_{0}} \Sigma E
\]
in which $\delta_{0}$ is a ghost.  We can then embed $\pi_{*}M_{1}$
into an injective module $I_{1}$, and get an analogous exact triangle
\[
M_{1} \xrightarrow{} J_{1} \xrightarrow{} M_{2} \xrightarrow{\delta_{1}} \Sigma M_{1}
\]
in which $\delta_{1}$ is a ghost.  Now the composite 
\[
M_{2} \xrightarrow{\delta_{1}} \Sigma M_{1} \xrightarrow{\delta_{0}}
\Sigma^{2} E
\]
is necessarily trivial, because $\rgldim E=1$.  It is represented in
the universal coefficient spectral sequence 
\[
E_{2}^{s,t} =\Ext_{E_{*}}^{s,t} (\pi_{*}M_{2},E_{*}) \Rightarrow
\cat{D} (E) (M_{2},E)_{t-s}
\]
on the $2$-line by the extension
\[
0 \xrightarrow{} R \xrightarrow{} I_{0} \xrightarrow{} I_{1}
\xrightarrow{} M_{2} \xrightarrow{} 0
\]
which is trivial if and only if $\idim R\leq 1$.  This element of
$E_{2}$ is a permanent cycle, since it represents a map, but it must
not survive the spectral sequence.  Therefore, there must be a
differential, necessarily a $d_{2}$, that hits it.  The source of such
a differential is a map $\pi_{*}M_{2}\xrightarrow{}E_{*}$.  But there
are no nonzero maps like this, because $M_{2}$ is the quotient of an
injective module.  Thus, the extension above must be $0$, so $\idim
R\leq 1$.  
\end{proof}

Now consider the case where $R$ is commutative Noetherian.  Then 
every injective module is a direct sum of copies of the
injective hulls $E (R/\ideal{p})$ of prime ideals $\ideal{p}$
(see~\cite[Section~3I]{lam}).  Every element of $E (R/\ideal{p})$ is
killed by $\ideal{p}^{n}$ for some $n$.  Thus, if $R$ is a domain, the
only possible injective module that can map to a projective module is
$E (R)$ itself.  But $E (R)$ is divisible, so every element can be
divided by every element of $R$.  Thus as long as $R$ is a Noetherian
domain that is not a field, there are no nonzero maps from an
injective to a projective module. 

We have therefore proved the following corollary. 

\begin{corollary}\label{cor-domain}
Suppose $E$ is an $S$-algebra with $\rgldim E=1$, and $E_{*}$ is a
commutative Noetherian domain that is not a field.  Then $E_{*}$ is
$1$-Gorenstein.
\end{corollary}


\providecommand{\bysame}{\leavevmode\hbox to3em{\hrulefill}\thinspace}
\providecommand{\MR}{\relax\ifhmode\unskip\space\fi MR }
\providecommand{\MRhref}[2]{%
  \href{http://www.ams.org/mathscinet-getitem?mr=#1}{#2}
}
\providecommand{\href}[2]{#2}

\end{document}